\documentclass[a4paper]{amsart}
\usepackage[margin=1.5in]{geometry}
\usepackage{amssymb,amsmath}
\address{Heng Xie, School of Mathematics, Sun-Yat Sen University, Guangzhou, China}
\email{xieh59@mail.sysu.edu.cn}
\usepackage{graphicx} 

\usepackage{booktabs}
\usepackage{multirow}
\usepackage{amsthm}
\usepackage{amsmath,amscd}
\usepackage{dsfont}
\usepackage{mathrsfs}
\usepackage[all,cmtip]{xy}
\usepackage{enumerate}
\usepackage[graphicx]{realboxes}
\usepackage{hyperref}
\theoremstyle{plain}
\newtheorem{theo}{Theorem}[section]
\usepackage{longtable}
\usepackage{lscape}
\usepackage{enumitem}
\usepackage{mathtools}
\usepackage{multirow}
\usepackage{floatrow}
\usepackage{tikz-cd}

\newcounter{mydiag}

\DeclareFontFamily{U}{mathx}{}
\DeclareFontShape{U}{mathx}{m}{n}{<-> mathx10}{}
\DeclareSymbolFont{mathx}{U}{mathx}{m}{n}
\DeclareMathAccent{\widehat}{0}{mathx}{"70}
\DeclareMathAccent{\widecheck}{0}{mathx}{"71}

\theoremstyle{theorem}
\newtheorem{lemma}[theo]{Lemma}
\newtheorem{proposition}[theo]{\bf Proposition}
\newtheorem{corollary}[theo]{\bf Corollary}

\theoremstyle{definition}
\newtheorem{remark}[theo]{Remark}

\newtheorem{definition}[theo]{Definition}
\newtheorem*{ackno}{Acknowledgement}

\newcommand{\h}{\mathscr{H}\kern -2pt om}
\newcommand{\id}{\textnormal{id}}
\newcommand{\op}{\textnormal{op}}
\newcommand{\tr}{\textnormal{tr}}
\newcommand{\Hom}{\textnormal{Hom}}
\newcommand{\quis}{\textnormal{quis}}
\newcommand{\hModh}{\textnormal{-Mod-}}

\newcommand{\hMod}{\textnormal{-Mod}}

\newcommand{\can}{\textnormal{can}}

\newcommand{\Fun}{\textnormal{Fun}}
\newcommand{\cone}{\textnormal{cone}}
\newcommand{\veet}{\flat}

\newcommand{\GW}{\mathrm{GW}}

\title{Hermitian $K$-theory of quadric hypersurfaces}

\author{Heng Xie}

\begin{document}
	\begin{abstract}
		Let $k$ be a commutative ring in which $2$ is invertible. We prove that the Hermitian $K$-theory of quadric hypersurfaces admits fibration sequences relating it to the base ring and to Clifford algebras equipped with various duality coefficients. All shifts and twists are taken into account.
	\end{abstract}
	\maketitle

	\section{Introduction}
	Let $(P,q)$ be a non-degenerate quadratic form of rank $n$ over a commutative ring $k$.\  Let $Q_d \subseteq \mathbb{P}(P)$ be the quadric hypersurface defined by the quadratic form $q$ where $d = n-2$. In the 1980s, Swan \cite{Swan} showed that the algebraic $K$-theory $K(Q_d)$ is isomorphic to the direct sum of $K(k)^d$ and $K(C_0(q))$. Here, $C_0(q)$ is the even Clifford algebra of $(P,q)$. 
	
	This paper aims to understand the Hermitian $K$-theory of the quadric hypersurface $Q_d$ over $k$. Unlike algebraic $K$-theory, Hermitian $K$-theory allows shifts coming from chain complexes and twists by line bundles.\ This added flexibility leads to a richer but more complicated theory, involving more cases to analyze. Assume that  $\sigma$ is the canonical involution on $C(q)$.\ The canonical involution on $C(q)$ restricts to an algebra involution on the even  part $C_0(q)$ of $C(q)$.\ 	Furthermore, the canonical involution $
	\sigma$ on $C(q)$ also restricts to an involution on the odd part $C_1(q)$ of $C(q)$, which makes $(C_1(q),\sigma)$ into a duality coefficient on $(C_0(q),\sigma)$ in the sense of \cite[Section 7.4]{S-Inven}.\  
	
	\begin{theo}\label{theo:main-theorem}
		Let $\mathscr{L}$ be a line bundle on $Q_d$.	The following statements hold in the stable homotopy category of spectra $\mathcal{SH}$. 
		\begin{enumerate}
			\item If $\mathscr{L} = \mathscr{O}(1-d)$, then there is a stable equivalence of spectra
			$$  \xymatrix{ K(k)^{\lfloor \frac{d}{2} \rfloor } \oplus \GW^{[i]}(\mathfrak{A}_{\mathscr{L}}) \ar[r]^-{\simeq} &  \GW^{[i]}(Q_d,\mathscr{L}). }  $$
			\item If $\mathscr{L} = \mathscr{O}(-d)$, then there is a homotopy fibration sequence of spectra
			$$  \xymatrix{ K(k)^{\lceil \frac{d}{2} \rceil -1} \oplus \GW^{[i]}(\mathfrak{A}_{\mathscr{L}})  \ar[r] &   \GW^{[i]}(Q_d,\mathscr{L})  \ar[r]^-{p_*} &  \GW^{[i-d]}(k, \Lambda^n P). 	}	$$
		\end{enumerate}
		Moreover, the spectrum $\GW^{[i]}(\mathfrak{A}_{\mathscr{L}})$ has the following properties:
		\begin{enumerate}
			\item [\textnormal{(a)}] If $[\mathscr{L}] = [\mathscr{O}(1)]$ in the picard group $\mathrm{Pic}(Q_d)/2$,  then we have a stable equivalence
			\[  \GW^{[i]}(\mathfrak{A}_{\mathscr{L}})  \cong \GW^{[i]}(C_0(q)_\sigma, C_1(q)_\sigma) \]
			which  is the Grothendieck-Witt spectrum of left projective $C_0(q)$-modules with the duality coefficient in $(C_1(q),\sigma)$.
			\item [\textnormal{(b)}] 	If $[\mathscr{L}] = [\mathscr{O}]$ in $\mathrm{Pic}(Q_d)/2$, then the spectrum $ \GW^{[i]}(\mathfrak{A}_{\mathscr{L}}) $ fits into a distinguished triangle
			\begin{equation}\label{eq:trace-seq}
				\xymatrix{   \GW^{[i]}(k)  \ar[r] &    \GW^{[i]}(\mathfrak{A}_{\mathscr{L}})   \ar[r] &   \GW^{[i+1]}(C_0(q)_\sigma)   \ar[r]^-{\eta \cup \tr }  & \GW^{[i]}(k)[1] }  
			\end{equation} 
			where $\eta \in \GW^{[-1]}_{-1}(k)$ is the Bott element and where $\tr$ is the trace map. 
		\end{enumerate} 
	\end{theo}

	To the best of my knowledge, Theorem~\ref{theo:main-theorem} has not appeared in the literature for the case $[\mathscr{L}] = [\mathscr{O}(1)]$, even in the setting of classical Witt groups. By taking negative homotopy groups of $\GW^{[i]}(Q_d, \mathscr{L})$ when $[\mathscr{L}] = [\mathscr{O}]$, Theorem~\ref{theo:main-theorem} recovers the corresponding statement for Witt groups (cf.\ \cite[Theorem~1.2]{Xie19}). In general, the distinguished triangle~\eqref{eq:trace-seq} does not split. 
	For instance, in the Witt group case, the trace map
	\[
	\tr \colon W(C_0(q)_\sigma) \longrightarrow W(k)
	\]
	is nonzero when $k=\mathbb{R}$ and $q$ is anisotropic; see \cite[Theorem~1.4]{Xie19}. 
	In contrast, if $Q$ is isotropic (that is, if it admits a rational point), the associated long exact sequence of  \eqref{eq:trace-seq} always splits and decomposes into short split exact sequences.
	
	\noindent\textbf{Convention}. Throughout the paper, we assume that $\frac{1}{2} \in k$. Explicit construction of symmetric forms will exploit this condition (cf. Appendix \ref{sec:symmform}). 
	
	\section{Setup}
	In this paper, we work in the framework of dg (resp. exact) categories with dualities (cf. \cite{S-Inven} and \cite{S-JPAA}), and notations will be taken from \cite{Xie19} unless redefined in the context. 
	\subsection{Dg $k$-modules}
	The purpose of this subsection is to fix notations on dg modules.
	\begin{definition}
		A dg $k$-module is a graded $k$-module $$M = \bigoplus_{m\in \mathbb{Z}} M^m $$ 
		together with a $k$-linear map $d: M \to M$, 
		called a differential such that $d(M^m) \subset M^{m+1}$ for all $i \in \mathbb{Z}$ and $d \circ d =0$. A morphism $f:M \to N$ of dg $k$-modules is a morphism of graded $k$-modules such that $f(M^m) \subset N^m$ and $f(dx) =d (f(x))$. 
	\end{definition}
	A dg $k$-module can also be presented as a cochain complex 
	$$ \cdots \to M^s \xrightarrow{d^s} M^{s+1} \to \cdots \to M^t \xrightarrow{d^t} M^{t+1} \cdots $$
	of $k$-modules. 
	
	\begin{definition} Let $M,N$ be dg $k$-modules. 
		\begin{enumerate}
			\item Define the tensor product 
			$$M \otimes N: = \bigoplus_{m \in \mathbb{Z}} (M \otimes N)^m  $$ by 
			$(M\otimes N)^m = \bigoplus_{p+q = m} M^p \otimes N^q$
			where the differential is given by $d(x\otimes y)= d x \otimes y + (-1)^{|x|} x \otimes dy$. 
			\item Define the function dg $k$-module (aka. internal hom) 
			$$[M,N] : = \bigoplus_{m \in \mathbb{Z}} [M,N]^m   $$ 
			by 
			$[M,N]^m = \prod_{-p+q = m}  \mathrm{Hom}_k(M^p,N^q)$  
			where the differential is given by $df = d_N \circ f -(-1)^{|f|}f \circ d_M $. 
		\end{enumerate}
	\end{definition}
	\begin{remark} Let $M $ be a dg $k$-module. 
		\begin{enumerate}
			\item Let $T$ be the dg $k$-module with $k$ in degree $-1$ and zero in all the other degrees with differential zero. We can form the shifted dg $k$-module $M[1]$ as the tensor product $T\otimes M$. With a simple calculation, we have $M[1]^m =  M^{m+1}$ with differential $d^{M[1]} = -d^M$. One can define $M[-1]$ similarly and define $M[i]$ iteratively for any $i \in \mathbb{Z}$. 
			\item  If $f:M\to N$ is a morphism of dg $k$-modules, then we can form the cone $C(f)$ of $f$ by setting $C(f) = N \oplus M[1] $ and the differential is defined by 
			\[N \oplus M[1] \xrightarrow{\begin{pmatrix}
					d^N & f \\
					0 & d^{M[1]}
			\end{pmatrix}} N \oplus M[1] \]
			Let $C$ be the dg $k$-module $k \langle 0 \rangle \oplus k \langle -1 \rangle$ with differential $d x = x$ for $x \in k$ with $|x| = -1$.	Note that $CM:= C\otimes M$ is precisely $C(M \xrightarrow{\mathrm{id}} M)$, and we have a cocartesian square
			\[\xymatrix{ M \ar[r] \ar[d]^-{f} & CM \ar[d] \\
				N \ar[r] & C(f).}\]
			\item Consider $k$ as a dg $k$-module concentrated in degree zero. Note that the dg $k$-module $[M, k[i]]$ has $[M, k[i]]^m = \Hom_k(M^{-i-m}, k) $, and the differential is given by $df = (-1)^{|f|+1}f \circ d^M$. It is not hard to check that $[M,k][i] = [M,k[i]]$. 
		\end{enumerate}
	\end{remark}
	
	\section{Clifford algebras and Tate resolution}\label{Cduality}
	\subsection{Quadric hypersurfaces}
	Let $(P,q)$ be a non-degenerate quadratic form of rank $n$ over $k$. Let $P^*$ denote the $k$-module $\Hom_k(P,k)$. Then, one can form the symmetric algebra $S(P^*)$.\ Let  $A=S(P^{*})/(q)$ be the homogeneous ring defined by $q$ viewing as an element in $S^2 (P^*)$, cf.\  \cite[Section 2]{Swan}.\ Define $Q$ to be the projective variety $\textnormal{Proj} A$, which is called the smooth quadric associated to $(P,q)$. By \cite{Swan}, $Q$ is smooth of relative dimension $d = n-2$ over $\textnormal{Spec\,}(k)$. 
	There is a pairing $\mathscr{O}(i) \otimes P^{*} \rightarrow \mathscr{O}(i+1)$ induced by the multiplication in the symmetric algebra.\ It follows that one has a map $ \mathscr{O}(i) \rightarrow \mathscr{O}(i+1) \otimes P   $ by the following composition
	$$  \mathscr{O}(i) \to  \Hom_k(P^{*}, \mathscr{O}(i+1))  \xrightarrow{\cong} \mathscr{O}(i+1) \otimes \Hom_k(P^{*}, k) \rightarrow  \mathscr{O}(i+1) \otimes P $$
	where the first map is induced by the pairing, the middle map is the canonical isomorphism, and the last map is induced by the double dual identification.
	
	\subsection{Clifford sequences}
	The Clifford algebra $ C(q)$ of the quadratic form $(P,q)$ is defined as $T(P)/I(q)$ where $T(P)$ is the tensor algebra of $P$ and where $I(q)$ is the two-sided ideal generated by $v \otimes v -q(v)$ for all $v \in P$.\ Moreover, the Clifford algebra $C(q) = C_0(q) \oplus C_1(q)$ has a $\mathbb{Z}/2\mathbb{Z}$-grading induced by the grading on the tensor algebra $T(P)$.\ 
	Let $r = d+1$.	If $M$ is a $\mathbb{Z}/2\mathbb{Z}$-graded \textit{left} $C(q)$-module, then we define a dg $\mathscr{O}$-module 
	$$\mathrm{Cliff}(M) = \bigoplus_{i \in \mathbb{Z}} \mathrm{Cliff} (M)^{-i} $$
	where 
	$$ \mathrm{Cliff} (M)^{-i} = \mathscr{O}(-i) \otimes M_{i+r} $$
	and the differential 
	$$d^{-i}: \mathrm{Cliff} (M)^{-i} \to \mathrm{Cliff} (M)^{-i+1}  $$
	is given by the composition 
	$$\mathscr{O}(-i) \otimes M_{i+r} \rightarrow \mathscr{O}(-i+1) \otimes P \otimes M_{i+r} \rightarrow \mathscr{O}(-i+1) \otimes M_{i+r-1}.$$
	If $M= C(q)$, then we note that the dg $\mathscr{O}$-module $\mathrm{Cliff}(C(q))$ 
	has a right $C_0(q)$-module structure. 
	
	\begin{definition}[Swan bundle]
		Define 
		$$ \mathscr{U}_{i}(M) := \textnormal{coker}\Big[ \mathrm{Cliff} (M)^{-i-2} \to \mathrm{Cliff} (M)^{-i-1}  \Big] $$
		Set  $\mathscr{U}_{i}:=\mathscr{U}_{i}(C(q))$.
	\end{definition}
	Note that $ \mathscr{U}_{i}$ has a right $C_0(q)$-module structure. Swan proved $$\textnormal{End}( \mathscr{U}_{i}) \cong C_0(q)$$ 
	for any $i \in \mathbb{Z}$, cf.\ \cite[Lemma 8.7]{Swan}, and $\mathscr{U}_i(M) = \mathscr{U}_{i}\otimes_{C_0(q)} M_0$.\ 
	
	\subsection{Tate resolution}
	Let $\Lambda^i = \Lambda^iP^*$ be the $i$-th exterior power and let $$\Lambda^{(l)} = \begin{cases}
		\Lambda^0 + \Lambda^2 + \Lambda^4 + \cdots& \textnormal{if $l$ is even} \\
		\Lambda^1 + \Lambda^3 + \Lambda^5 + \cdots  & \textnormal{if $l$ is odd} \\
	\end{cases}$$ 
	Define $$\mathbf{C}_{ } :=\bigoplus_{i \in \mathbb{Z}} \mathbf{C}^{-i} $$
	where 
	$$\mathbf{C}^{-i} = \mathscr{O}(-i) \otimes \Lambda^{(i+1)}. $$
	The differential 
	$ \Delta^{-i} : \mathbf{C}^{-i} \to \mathbf{C}^{-i+1} $
	is given by $\Delta_1^{-i} + \Delta_2^{-i}$ where 
	$$\Delta^{-i}_1 (f \otimes p_1 \wedge \cdots \wedge p_k ) = \sum_s (-1)^{s+1} f p_s \otimes (p_1 \wedge \cdots \wedge \hat{p}_s \wedge \cdots \wedge p_k)$$
	and 
	$$\Delta^{-i}_2(f \otimes \omega ) = \sum_i f \xi_i \otimes (\beta_i \wedge \omega)$$
	Here, $\xi_i, \beta_i \in P^*$ are finite elements such that $\gamma=\sum_i \xi_i \otimes \beta_i \in P^* \otimes P^*$ maps to $q \in S_2 (P^*)$ via the natural surjection $P^* \otimes P^* \to S_2(P^*)$. Note that $\Lambda := \oplus_{m \in \mathbb{Z}} \Lambda^m= \Lambda^{(0)} \oplus \Lambda^{(1)}$ is a $\mathbb{Z}/2 \mathbb{Z}$-graded $C(q)$-algebra. Define the $\mathbb{Z}/2$-graded $C(q)$-algebra $\Lambda(d)$ as by the $d$-th shift of $\Lambda$, i.e. $\Lambda(d)_{i} = \Lambda^{(d+i)} $. 
	\begin{lemma}
		$(\mathbf{C}_{ }, \Delta) \cong \mathrm{Cliff}(\Lambda(d)).$
	\end{lemma}
	\begin{proof}
 See \cite[Lemma 8.4, p.\ 131]{Swan}.
	\end{proof}
	Let $$\Lambda^{\langle l \rangle} = \begin{cases}
		\Lambda^l + \Lambda^{l-2} + \cdots + \Lambda^2 + \Lambda^0 & \textnormal{if $l$ is even} \\
		\Lambda^l + \Lambda^{l-2} + \cdots + \Lambda^3 + \Lambda^1 & \textnormal{if $l$ is odd.}
	\end{cases}$$ 
	These settings yield a short split exact sequence
	\begin{equation}
		\xymatrix{ 0 \ar[r] & \Lambda^{\{l\}}\ar[r]_-{\iota} &  \Lambda^{(l)} \ar[r]_-{\pi} \ar@/_.8pc/[l]_-{\rho} & \ar@/_.8pc/[l]_{\tau} \Lambda^{\langle l \rangle} \ar[r]  & 0 }
	\end{equation}
	with $\rho\tau = 0$,	where 
	$$\Lambda^{\{l\}} = 		\Lambda^{l+2} +  \Lambda^{l+4}  +  \Lambda^{l+6}+ \cdots
	$$ 
	\begin{definition}
		Define  
		$$\mathbf{T}_{ } := \bigoplus_{i \geq -1} \mathbf{T}^{-i}, \quad \quad \widetilde{\mathbf{T}}_{ }:=  \bigoplus_{i \leq d-1} \widetilde{\mathbf{T}}^{-i}  $$
		where 
		$$\mathbf{T}^{-i} = \mathscr{O}(-i) \otimes \Lambda^{\langle i+1 \rangle}, \quad \quad \widetilde{\mathbf{T}}^{-i} = \mathscr{O}(-i) \otimes \Lambda^{\{i+1\}} $$
		The differential 
		$ \partial^{-i} : \mathbf{T}^{-i} \to \mathbf{T}^{-i+1} $
		is given by the composition 
		$$ \mathscr{O}(-i) \otimes \Lambda^{\langle i+1 \rangle} \xrightarrow{1 \otimes \tau} \mathscr{O}(-i) \otimes \Lambda^{(i+1)}  \xrightarrow{\Delta} \mathscr{O}(-i+1) \otimes \Lambda^{(i)}  \xrightarrow{1 \otimes \pi} \mathscr{O}(-i+1) \otimes \Lambda^{\langle i \rangle},  $$
		and 
		the differential 
		$ \widetilde{\partial}^{-i} : \widetilde{\mathbf{T}}^{-i} \to \widetilde{\mathbf{T}}^{-i+1} $
		is given by the composition 
		$$ \mathscr{O}(-i) \otimes \Lambda^{\{ i+1 \}} \xrightarrow{1 \otimes \iota} \mathscr{O}(-i) \otimes \Lambda^{(i+1)}  \xrightarrow{\Delta} \mathscr{O}(-i+1) \otimes \Lambda^{(i)}  \xrightarrow{1 \otimes \rho} \mathscr{O}(-i+1) \otimes \Lambda^{\{ i \}}.  $$
	\end{definition}    
	\begin{lemma}
		The modules $(\mathbf{T}_{ }, \partial)$ and $(\widetilde{\mathbf{T}}, \widetilde{\partial})$ are dg $\mathscr{O}$-modules, and the sequence
		\begin{equation}\label{eq:dgomoduleT}
			\xymatrix{ 0 \ar[r] & \widetilde{\mathbf{T}}_{ }\ar[r]^-{\iota} & \mathbf{C}_{ } \ar[r]^-{\pi} &  \mathbf{T}_{ }  \ar[r] & 0 }
		\end{equation} 
		is an exact sequence of dg $\mathscr{O}$-module that splits degreewise.\
	\end{lemma}
	\begin{proof}
		The result follows directly from the commutativity of the following diagram
		\begin{equation*}
			\xymatrix{
				\mathbf{C}^{-i} \ar[r]^-{\Delta} \ar[d]_{\cong}^-{\left(\begin{smallmatrix}
						\pi  \\
						\rho 
					\end{smallmatrix}\right)} & \mathbf{C}^{-i+1} \ar[d]_{\cong}^-{\left(\begin{smallmatrix}
						\pi  \\
						\rho 
					\end{smallmatrix}\right)} \\
				\mathbf{T}^{-i} \oplus \widetilde{\mathbf{T}}^{-i} \ar[r]^-{\left(\begin{smallmatrix}
						\partial & 0 \\
						\zeta & \widetilde{\partial}
					\end{smallmatrix}\right)} & \mathbf{T}^{-i+1} \oplus  \widetilde{\mathbf{T}}^{-i+1} }
		\end{equation*}
		where $\zeta = \rho \Delta \tau$. Note that, inside $\Lambda^{(l)}$, $\Lambda^{\{l\}}$ or $\Lambda^{\langle l\rangle}$, each different summand $\Lambda^i$ and $\Lambda^j$ must satisfy $|i-j| \geq 2$. However, $\Delta= \Delta_1 + \Delta_2$ only increase or decrease the degree by $\pm 1$. This shows that there is no map from $\widetilde{\mathbf{T}}^{-i}$ to  $\mathbf{T}^{-i+1}$, i.e. $\pi\Delta\iota =0$. Since $\Delta \circ \Delta =0$, we see that $\partial \circ \partial =0$ and $\widetilde{\partial}\circ \widetilde{\partial} = 0$. 
	\end{proof} 
	\begin{remark}
		The dg $\mathcal{O}$-module $\mathbf{T}_{ }$ is precisely the Tate resolution given in \cite{Swan}. 
	\end{remark}
	
	\begin{remark}We have  $\mathbf{T}^1= \mathscr{O}(-1), \widetilde{\mathbf{T}}^{1-d}= \mathscr{O}(1-d)\otimes \Lambda^{n} $, since $\Lambda^{\{d\}}= \Lambda^n, \Lambda^{\langle 0 \rangle} =k$.
	\end{remark}
	\subsection{Modification}	Let $\mathscr{H}:=\mathscr{U}_{d-2}(\Lambda(d)) $ and $\mathscr{G}:= \mathscr{U}_{-1}(\Lambda(d))$. By definition, the dg $\mathscr{O}$-module $\mathbf{C}_{ }$ can be truncated into three acyclic chain complexes
	\[\cdots \to  \mathbf{C}^{-d} \to \mathbf{C}^{1-d} \xrightarrow{\beta} \mathscr{H}   \to 0, \]
	\[0 \to \mathscr{H} \xrightarrow{\alpha} \mathbf{C}^{2-d} \to \mathbf{C}^{3-d} \to \cdots \to \mathbf{C}^{-1} \to \mathbf{C}^{0} \xrightarrow{\gamma} \mathscr{G} \to 0,   \]
	\[0 \to \mathscr{G} \xrightarrow{\eta} \mathbf{C}^{1} \to \mathbf{C}^{2} \to \cdots \]
	It follows from \cite[Section 4]{Xie19} that  the two bounded chain complexes
	\[\mathbf{F } = (0 \to \widetilde{\mathbf{T}}^{1-d} \xrightarrow{\beta\iota} \mathscr{H} \xrightarrow{\pi\alpha} \mathbf{T}^{2-d} \xrightarrow{\partial} \cdots \xrightarrow{\partial} \mathbf{T}^{-1} \xrightarrow{\partial} \mathbf{T}^{0} \xrightarrow{\partial} \mathbf{T}^{1} \to 0)\]
	\[\widetilde{\mathbf{F }} = (0 \to \widetilde{\mathbf{T}}^{1-d} \xrightarrow{\widetilde{\partial}}  \widetilde{\mathbf{T}}^{2-d} \xrightarrow{\widetilde{\partial}}  \widetilde{\mathbf{T}}^{3-d} \xrightarrow{\widetilde{\partial}} \cdots \xrightarrow{\widetilde{\partial}}  \widetilde{\mathbf{T}}^{0} \xrightarrow{\gamma\iota} \mathscr{G} \xrightarrow{\pi\eta} \mathbf{T}^{1} \to 0)\]
	are acyclic. It is helpful to summarize the sequence (\ref{eq:dgomoduleT}) and our notations by the following ladder diagrams of chain complexes.
	\begin{equation}\label{eq:explicit-TCT}
		\xymatrix@C=20pt{ 
			0 \ar[r] \ar[d] & \widetilde{\mathbf{T}}^{1-d} \ar[dr]^-{\beta\iota} \ar[rr]^-{\tilde{\partial}} \ar[d]^-{\iota} && \widetilde{\mathbf{T}}^{2-d} \ar[r]^-{\tilde{\partial}} \ar[d]^-{\iota} &\cdots \ar[r]^-{\tilde{\partial}} & \widetilde{\mathbf{T}}^{-1} \ar[r]^-{\tilde{\partial}} \ar[d]^-{\iota}& \widetilde{\mathbf{T}}^{0} \ar[dr]^-{\gamma \iota} \ar[rr]^-{\tilde{\partial}} \ar[d]^-{\iota} && \widetilde{\mathbf{T}}^{1} \ar[r]^-{\tilde{\partial}} \ar[d]^-{\iota} & \widetilde{\mathbf{T}}^{2} \ar[d]^-{=} \\  
			\mathbf{C}^{-d} \ar[r]^-{\Delta} \ar[d]^-{=} & \mathbf{C}^{1-d} \ar@{->>}[r]^-{\beta} \ar[d]^-{\pi}  &\ar@{>->}[r]^-{\alpha} \mathscr{H} \ar[dr]_-{\pi\alpha} & \mathbf{C}^{2-d} \ar[r]^-{\Delta} \ar[d]^-{\pi} &\cdots \ar[r]^-{\Delta}& \mathbf{C}^{-1} \ar[r]^-{\Delta} \ar[d]^-{\pi} & \mathbf{C}^{0} \ar[r]^-{\gamma} \ar[d]^-{\pi} & \ar[r]^-{\eta} \ar[dr]_-{\pi \eta} \mathscr{G}  & \mathbf{C}^{1} \ar[r]^-{\Delta} \ar[d]^-{\pi} & \mathbf{C}^{2} \ar[d] \\  
			\mathbf{T}^{-d} \ar[r]^-{\partial} & \mathbf{T}^{1-d} \ar[rr]^-{\partial} && \mathbf{T}^{2-d} \ar[r]^-{\partial} &\cdots \ar[r]^-{\partial} & \mathbf{T}^{-1} \ar[r]^-{\partial} & \mathbf{T}^{0}  \ar[rr]^-{\partial} && \mathbf{T}^{1} \ar[r] & 0
	}\end{equation}

	\section{Construction of symmetric forms}	
	\subsection{Symmetric forms on Tate resolution}\label{sub:tate-dual}
	Let $\mathscr{I}:= \mathcal{O}(2-d) \otimes \Lambda^nP^*$. Consider the category with duality $(\mathscr{E}, *_{\mathscr{I}},\eta_{\mathscr{I}})$ where $\mathscr{E} =\mathrm{Vect}(Q)$ is the category 
	of finite rank vector bundles over $Q$, and 
	$*_{\mathscr{I}}$ is the duality $\h_\mathscr{O}(-,\mathscr{I}) $ with the canonical double dual identification $\eta_{\mathscr{I}}$.\ The morphism of $k$-modules
	\begin{align*}
		\Lambda^{i}P^* \otimes \Lambda^{n-i} P^* & \to  \Lambda^nP^* \\
		x \otimes y & \mapsto x \wedge y
	\end{align*}
	induces a morphism of $\mathcal{O}$-modules
	\[ \psi: \mathbf{C}^{-i} \to (\mathbf{C}^{2-d+i})^{*_{\mathscr{I}}} \]
	by setting  
	\[\psi(f\otimes x)(g \otimes y) =
	\begin{cases}
		(-1)^{\frac{s(s-1)}{2}}fg \otimes x \wedge y & \textnormal{if $s+t =n$ } \\
		0 & \textnormal{if otherwise }
	\end{cases} \]
	for all $f \in \mathcal{O}(-i),g\in \mathcal{O}(2-d+i), x= x_1 \wedge \cdots \wedge x_s \in \Lambda^{(i+1)}, y = y_1 \wedge \cdots \wedge y_t \in \Lambda^{(d-1-i)}$. 
	
	\begin{lemma}\label{lem:psi}
		The morphism $\psi$ induces a symmetric space 
		\[\xymatrix{ \mathbf{C} \ar[d]^-{\psi} & \cdots \ar[r] &  \mathbf{C}^{-i}  \ar[d]^-{\psi} \ar[r]^-{\Delta} &\cdots \ar[r]^-{\Delta}  & \mathbf{C}^{-1} \ar[d]^-{\psi}  \ar[r]^-{\Delta} & \mathbf{C}^{0} \ar[d]^-{\psi} \ar[r] & \cdots  \\
			\quad\,\, \mathbf{C}^{*_{\mathscr{I}}^{d-2}} & \cdots \ar[r] & (\mathbf{C}^{2-d+i})^{*_{\mathscr{I}}} \ar[r]^-{\Delta^*} &\cdots \ar[r]^-{\Delta^*} & (\mathbf{C}^{3-d})^{*_{\mathscr{I}}}  \ar[r]^-{\Delta^*} & (\mathbf{C}^{2-d})^{*_{\mathscr{I}}} \ar[r] & \cdots
		} \]
		in the category with duality $(\mathscr{O}\hMod, *_{\mathscr{I}}^{d-2},\eta_{\mathscr{I}}^{d-2})$.
	\end{lemma}
	\begin{proof}
		Let us check the commutativity of the diagram
		\begin{equation}\label{eq:C-chaincomplex} \xymatrix{ \mathbf{C}^{-i} \ar[r]^-{\Delta} \ar[d]^-{\psi} & \mathbf{C}^{-i+1} \ar[d]^-{\psi} \\
				(\mathbf{C}^{2-d+i})^{*_{\mathscr{I}}} \ar[r]^-{\Delta^*} & (\mathbf{C}^{3-d+i})^{*_{\mathscr{I}}}.
		} \end{equation}
		This amounts to showing that the equality $$\psi \Delta (f \otimes x) (g\otimes z) = \Delta^*\psi (f \otimes x) (g\otimes z)$$
		holds for any $x = x_1 \wedge \dots \wedge x_s \in \Lambda^{(i+1)}$ and $ z = z_1 \wedge \dots \wedge z_r \in \Lambda^{(d-1-i)}$. 
		Define
		\[
		\hat{x}_k = x_1 \wedge \dots \wedge \hat{x}_k \wedge \dots \wedge x_s, \quad \hat{z}_j = z_1 \wedge \dots \wedge \hat{z}_j \wedge \dots \wedge z_r,
		\]
		with $x_j$ (resp. $z_k$) deleted from $x$ (resp. $z$). Recall that $\Delta = \Delta_1 + \Delta_2$. Moreover, 
		\[
		\begin{array}{rll}
			\psi \Delta_1 (f \otimes x) (g\otimes z)& = \psi(\sum_{k=1}^s (-1)^{k+1} f x_k \otimes \hat{x}_k)(g\otimes z) \\
			&= (-1)^{\frac{(s-1)(s-2)}{2}}\sum_{k=1}^{s}(-1)^{k+1}fgx_k\otimes \hat{x}_k \wedge z \\
			&= (-1)^{\frac{(s-1)(s-2)}{2}}(-1)^{s-1}\sum_{j=1}^{r}(-1)^{j+1}fgz_j\otimes x \wedge \hat{z}_j \\
			&= (-1)^{\frac{(s-1)s}{2}}\sum_{j=1}^{r}(-1)^{j+1}fgz_j\otimes x \wedge \hat{z}_j \\
			& =\psi(f\otimes x)(\sum_{j=1}^{r}(-1)^{j+1}gz_j\otimes \hat{z}_j) \\
			& = \Delta_1^*\psi (f \otimes x) (g\otimes z)
		\end{array}\]
		for $s+r = n+1$, where the third equality holds since $\Delta_1(x\wedge z) =0$. 
		\[
		\begin{array}{rll}
			\psi \Delta_2 (f \otimes x) (g\otimes z)& = \psi(\sum_{t} f \xi_t \otimes (\beta_t \wedge x))(g\otimes z) \\
			& = (-1)^{\frac{s(s+1)}{2}} \sum_{t} f g \xi_t \otimes (\beta_t \wedge x \wedge z) \\
			& = (-1)^{\frac{s(s+1)}{2}}(-1)^{s} \sum_{t} f g \xi_t \otimes (x \wedge \beta_t \wedge z) \\
			&= (-1)^{\frac{s(s-1)}{2}}\sum_{t} f g \xi_t \otimes (x \wedge \beta_t \wedge z) \\
			&= \psi(f\otimes x) (\Delta_2(g \otimes z)) \\
			&=\Delta_2^* \psi(f\otimes x)(g \otimes z)
		\end{array}\]
		for $s+r = n-1$. This shows that Diagram \eqref{eq:C-chaincomplex} commutes. Finally, we check the symmetry. For  $x=x_1 \wedge \cdots \wedge x_s \in \Lambda^{(i+1)}, y = y_1 \wedge \cdots \wedge y_t \in \Lambda^{(d-1-i)}$
		\[\begin{array}{rll}
			\psi(f \otimes x) (g\otimes y)& = (-1)^{\frac{s(s-1)}{2}} fg \otimes x\wedge y \\
			&=(-1)^{\frac{s(s-1)}{2}}(-1)^{st} fg \otimes y\wedge x \\
			&=(-1)^{\frac{n(n-1)}{2}}(-1)^{\frac{t(t-1)}{2}}fg\otimes y\wedge x \\
			&=(-1)^{\frac{n(n-1)}{2}}\psi(g\otimes y)(f\otimes x) \\
			&=(-1)^{\frac{(d-2)(d-3)}{2}}\psi(g\otimes y)(f\otimes x)
		\end{array}\]
		whenever $s+t =n$. The result follows after comparing the sign conventions. 	\end{proof}
	
	\begin{lemma}\label{lem:phi}
		The map of complexes
		\[\xymatrix@C=1.3pc{ \mathbf{F } \ar[d]^-{\phi} & \widetilde{\mathbf{T}}^{1-d} \ar[d]^-{\phi_{-d}} \ar[r]^-{\beta \iota} & \mathscr{H} \ar[d]^-{\phi_{1-d}} \ar[r]^-{\pi \alpha} & \mathbf{T}^{2-d} \ar[d]^-{\phi_{2-d}} \ar[r]^-{\partial} &\cdots \ar[r] & \mathbf{T}^{-1} \ar[d]^-{\phi_{-1}} \ar[r]^-{\partial}  & \mathbf{T}^{0} \ar[d]^-{\phi_0}  \ar[r]^-{\partial} & \mathbf{T}^{1} \ar[d]^-{\phi_1}   \\
			\quad\,\, \mathbf{F }^{*_{\mathscr{I}}^{d-1}} & ({\mathbf{T}}^{1})^{*_{\mathscr{I}}} \ar[r]^-{\partial^*} & ({\mathbf{T}}^{0})^{*_{\mathscr{I}}} \ar[r]^-{\partial^*} & ({\mathbf{T}}^{-1})^{*_{\mathscr{I}}} \ar[r]^-{\partial^*} &\cdots \ar[r]^-{\partial^*} &(\mathbf{T}^{2-d})^{*_{\mathscr{I}}} \ar[r]^-{\alpha^*\pi^*} & {\mathscr{H}}^{*_{\mathscr{I}}}  \ar[r]^-{\iota^* \beta^*} & (\widetilde{\mathbf{T}}^{1-d})^{*_{\mathscr{I}}} 
		} \]
		is a (-1)-symmetric form in the category with duality $(\mathrm{Ch}^b(\mathcal{E}), *_{\mathscr{I}}^{d-1},\eta_{\mathscr{I}}^{d-1})$, where
		\[\phi_i := \begin{cases}
			(-1)^{1-d}\tau^*\psi\iota & \textnormal{for $i = -d$} \\ 
			(-1)^{1-d}\tau^* \psi\alpha & \textnormal{for $i = 1-d$} \\
			(-1)^{i}\tau^*\Delta^*\psi\tau & \textnormal{for $ 2-d \leq i \leq -1$} \\
			\alpha^* \psi \tau & \textnormal{for $i = 0$} \\
			\iota^*\psi\tau & \textnormal{for $i = 1$} 
		\end{cases} \]
		Moreover, the cone $C(\phi)$ is an acyclic chain complex. 
	\end{lemma}
	\begin{proof}
		Let us check that the map $\phi:  \mathbf{F} \to \mathbf{F}^{*_{\mathscr{I}}^{d-1}}$ is a morphism of chain complexes.	The commutativity of each square of $\phi$ follows from the following identities:
		\begin{align}
			\partial^*	\phi_{-d}& = \phi_{1-d} \beta \iota \label{6} \\
			\partial^* \phi_{1-d}&= \phi_{2-d} \pi\alpha  \label{7}  \\
			\partial^* \phi_i & = \phi_{i+1} \partial 
			\label{8} \\
			\phi_0 \partial & = \alpha^* \pi^* \phi_{-1} 
			\label{9} \\
			\phi_1 \partial & = \iota^* \beta^* \phi_0 
			\label{10}
		\end{align} 
		for any  $2-d \leq  i \leq -2$. 
		To check these identities, it is helpful to draw the dual of the diagram \eqref{eq:explicit-TCT} as follows:
		\[
		\xymatrix@C=20pt{ 
			({\mathbf{T}}^{1})^{*_{\mathscr{I}}}\ar[dr]^-{\eta^*\pi^*} \ar[rr]^-{{\partial}^*} \ar[d]^-{\pi^*} && ({\mathbf{T}}^{0})^{*_{\mathscr{I}}}\ar[r]^-{{\partial}^*} \ar[d]^-{\pi^*} &\cdots \ar[r]^-{{\partial}^*} &  ({\mathbf{T}}^{2-d})^{*_{\mathscr{I}}}\ar[dr]^-{\alpha^*\pi^*} \ar[rr]^-{{\partial}^*} \ar[d]^-{\pi^*} && ({\mathbf{T}}^{1-d})^{*_{\mathscr{I}}} \ar[d]^-{\pi^*}  \\  
			(\mathbf{C}^{1})^{*_{\mathscr{I}}}\ar@{->>}[r]^-{\eta^*} \ar[d]^-{\iota^*}  &\ar@{>->}[r]^-{\gamma^*} \mathscr{G}^{*_{\mathscr{I}}}\ar[dr]_-{\iota^*\gamma^*} & (\mathbf{C}^{0})^{*_{\mathscr{I}}}\ar[r]^-{\Delta^*} \ar[d]^-{\iota^*} &\cdots \ar[r]^-{\Delta^*}&  (\mathbf{C}^{2-d})^{*_{\mathscr{I}}}\ar[r]^-{\alpha^*}\ar[d]^-{\iota^*} & \ar[r]^-{\beta^*} \ar[dr]_-{\iota^*\beta^*} {\mathscr{H}}^{*_{\mathscr{I}}} & (\mathbf{C}^{1-d})^{*_{\mathscr{I}}}  \ar[d]^-{\iota^*}  \\  
			(\widetilde{\mathbf{T}}^{1})^{*_{\mathscr{I}}}\ar[rr]^-{\tilde{\partial}^*} && (\widetilde{\mathbf{T}}^{0})^{*_{\mathscr{I}}}\ar[r]^-{\tilde{\partial}^*} &\cdots \ar[r]^-{\tilde{\partial}^*} &  (\widetilde{\mathbf{T}}^{2-d})^{*_{\mathscr{I}}} \ar[rr]^-{\tilde{\partial}^*} && (\widetilde{\mathbf{T}}^{1-d})^{*_{\mathscr{I}}}
		}\]  
		The following formulas play important roles in checking Formula \eqref{6} - \eqref{10}. 
		\begin{align}
			\psi \tau \pi  & = \iota \rho \psi \label{11} \\
			\iota \rho + \tau \pi &=1 \label{12}  \\
			\pi^* \tau^* \psi  & = \psi \iota \rho
			\label{13} \\
			\iota^* \psi \tau \pi &= \iota^* \psi \label{14} \\
			\pi^*\tau^*\psi\iota & = \psi\iota  \label{15}\\
			\alpha \beta &= \Delta  
			\label{16} \\
			\psi \Delta  & = \Delta^* \psi 
			\label{17}
		\end{align} 
		where Formulas \eqref{11} - \eqref{15} are induced by the construction of the split exact sequence \eqref{eq:dgomoduleT} and the definition of $\psi$.  Formula \eqref{16} follows from the definition of $\alpha$ and $\beta$ and Formula \eqref{17} is deduced by Lemma \ref{lem:psi}.
		%
		%
		Let us check Formula	\eqref{8} by the following equalities
		\[			
		\begin{array}{rll}
			\partial^*	\phi_{i}& =  (-1)^{i} \tau^* \Delta^* \pi^* \tau^* \Delta^* \psi \tau    \\
			& =  (-1)^{i} \tau^* \Delta^* \pi^* \tau^* \psi \Delta  \tau    \\
			& =  (-1)^{i} \tau^* \Delta^*  \psi \iota \rho \Delta \tau   \\
			& =  (-1)^{i} \tau^* \Delta^*  \psi (1 - \tau \pi) \Delta \tau   \\
			& =  (-1)^{i+1} \tau^* \Delta^*  \psi   \tau \pi \Delta \tau   \\
			&= \phi_{i+1} \partial
		\end{array}	 \]			
		Formulas \eqref{6}, \eqref{7}, \eqref{9}, and \eqref{10} can be  obtained in a similar way. 
		Next, we check the symmetry: 
		\[ \phi =  - \phi^{*^{d-1}} \eta^{d-1} \]
		For $ 2-d \leq i \leq -1$, we have
		\[	\begin{array}{rll}
			-	(\phi^{*^{d-1}})_i (\eta^{d-1})_i & =- (-1)^{\frac{(d-1)(d-2)}{2}} \phi_{-i-(d-1)}^{*} \eta_T  \\
			& = - (-1)^{\frac{(d-1)(d-2)}{2}} \phi_{1-d-i}^{*} \eta_T \\
			& =- (-1)^{\frac{(d-1)(d-2)}{2}} (-1)^{1-d-i} \tau^*  \psi^* \Delta^{**} \tau^{**} \eta_T \\
			& =- (-1)^{\frac{(d-1)(d-2)}{2}} (-1)^{1-d-i} \tau^*  \psi^* \eta_C \Delta \tau\\
			& =- (-1)^{\frac{(d-1)(d-2)}{2}} (-1)^{1-d-i} (-1)^{\frac{(d-2)(d-3)}{2}} \tau^*  \psi \Delta \tau\\
			& = (-1)^{i} \tau^*  \Delta^* \psi \tau\\
			& = \phi_i
		\end{array} \]
		The symmetry for $i=0,1$ can be obtained similarly.\	This shows that the form $\phi$ is a (-1)-symmetric form in the category with duality $(\mathrm{Ch}^b(\mathscr{E}), *^{d-1}, \eta^{d-1})$.
	\end{proof}    
	
	Let $\mathscr{K}:= \mathscr{O}(-d)\otimes \Lambda^n P^*$. Next, we wish to tensor the chain complex $\mathbf{F}$ with $\mathscr{O}(-1)$ so that it starts with $\mathscr{O}$ concentrated in degree $1$ and respects the duality $*_\mathscr{K} = [-, \mathscr{K}]_{\mathscr{O}}$. Let $\mathbf{E} := \mathbf{F}(-1)$. Applying Lemma \ref{lemma:symm-E} to $\mathbf{E}$, we immediately get the following result. 
	\begin{corollary}
		\label{prop:Kphi}
		There is a symmetric form
		$$ \xymatrix@C=15pt{
			\mathbf{K} \ar[d]_-{{\varphi }} &  \mathbf{T}^{1-d}(-1) \ar[d]_-{\varphi_{-d}} \ar[r]  &  \mathscr{H}(-1) \ar[d]_-{\varphi_{-n+1}}  \ar[r] &   \cdots \ar[r] &  \mathbf{T}^{-1}(-1) \ar[d]_-{\varphi_{-1}} \ar[r]^-{\partial}  & \mathbf{T}^{0}(-1) \ar[d]_-{\varphi_{0}} \\
			\quad	\mathbf{K}^{*_{\mathscr{K}}^{n}} &  \mathbf{T}^{0}(-1)^{*_{\mathscr{K}}} \ar[r]^-{\partial^{*}}  &  \mathbf{T}^{-1}(-1)^{*_{\mathscr{K}}} \ar[r] &   \cdots \ar[r] &  \mathscr{H}(-1)^{*_{\mathscr{K}}} \ar[r]  & \mathbf{T}^{1-d}(-1)^{*_{\mathscr{K}}} }$$ 
		in the category with duality $(\mathrm{Ch}^b(\mathcal{E}), *_{\mathscr{K}}^d, \eta_{\mathscr{K}}^d)$. 
	\end{corollary}

	\begin{remark}
		The symmetry of the form ${\varphi }$ could also be checked by a less direct argument provided in \cite[Proposition 7.5]{Xie19}, i.e. reducing the problem to the base field $\mathbb{C}$ and compare to the computation made in \cite{Zibrowius}. 
	\end{remark}
	
	\subsection{Trivial duality} 
	Let $\ell = 2a$ and let $\mathscr{L}= \mathscr{O}(-\ell)$.	Consider the category with duality
	\[(\mathscr{O} \hModh C_0(q), \sharp_\mathscr{L}, \can_\mathscr{L}) \]
	where $\sharp_\mathscr{L}(M) = [M^\op, \mathscr{L}]_\mathscr{O}$ for any $M \in \mathscr{O} \hModh C_0(q)$, and where $\can_\mathscr{L}: M \to M^{\sharp_\mathscr{L} \sharp_\mathscr{L}} $ is given by	
	$\can(x)(f^\op) = (-1)^{|f||x|} i(f(x^\op))$. 
	
	Note that we have an isomorphism of exact sequences
	\begin{equation}\label{eq:trivial-duality}
		\xymatrix{
			\cdots \ar[r]   & \mathscr{O}(-i-2)\otimes C_{i+r+2}(q) \ar[d]^{\eta} \ar[r] & \mathscr{O}(-i-1)\otimes C_{i+r+1}(q) \ar[d]^{\eta} \ar[r] & \cdots \\
			\cdots \ar[r]  & (\mathscr{O}(i+2-\ell) \otimes C_{i+r+2}(q))^{\sharp_\mathscr{L}} \ar[r] & (\mathscr{O}(i+1-\ell) \otimes C_{i+r+1}(q))^{\sharp_\mathscr{L}}  \ar[r] & \cdots
		}
	\end{equation}
	of $\mathscr{O}$-modules and right $C_0(q)$-modules, where we define
	\[ \eta( f \otimes x)(g \otimes y^\op) = fg~ \tr(\sigma(x) \cdot y) . \]
	By taking cokernels of horizontal maps in \eqref{eq:trivial-duality}, we get isomorphisms
	\[ h_i: \mathscr{U}_{i} \longrightarrow (\mathscr{U}_{-i-1+\ell})^{\sharp_\mathscr{L}}  . \]
	\begin{lemma}\label{U0U1dual}  Let $\ell = 2a$. Then, there is a symmetric space of exact sequences 
		$$\xymatrix{ \mathscr{U}_{a} \ar[r] \ar[d]^-{h_{a}} &  \mathscr{O}(-a) \otimes C_{r+a}(q) \ar[r] \ar[d]^-{\eta} & \mathscr{U}_{a-1} \ar[d]^-{h_{a-1}} \\
			(\mathscr{U}_{a-1})^{\sharp_\mathscr{L}} \ar[r] &  (\mathscr{O}(-a) \otimes C_{r+a}(q))^{\sharp_\mathscr{L}} \ar[r]^-{m^\vee} & (\mathscr{U}_{a})^{\sharp_\mathscr{L}} }$$
		of  $\mathscr{O}$-modules and right $C_0(q)$-modules.\ \end{lemma}
	\begin{proof}
		See \cite[Lemma 5.12]{Xie19} for the case $a=0$. For any $m$, the proof is the same. 
	\end{proof}
	\begin{proposition}\label{prop:clmu}
		There is a symmetric form
		\begin{equation}
			\xymatrix{
				Cl \ar[d]_-{\mu} & \mathscr{U}_{a} \ar[r] \ar[d] & \mathscr{O}(-a) \otimes C_{r+a}(q) \ar[d] \\
				\quad Cl^{\sharp^{-1}_\mathscr{L}}  & (\mathscr{O}(-a) \otimes C_{r+a}(q))^{\sharp_\mathscr{L}} \ar[r] & (\mathscr{U}_{a})^{\sharp_\mathscr{L}}	}	
		\end{equation}
		in the category with duality $(\mathscr{O} \hModh C_0(q), \sharp^{-1}_\mathscr{L},\can_{\sharp_\mathscr{L}}^{-1}) $.
	\end{proposition}
	\begin{proof}
		See \cite[Lemma 6.3]{Xie19}. Alternatively, one can use the equivalence 
		\[(\mathscr{O} \hModh C_0(q), \sharp^{-1}_\mathscr{L},\can_{\sharp_\mathscr{L}}^{-1}) \cong (\mathscr{O} \hModh C_0(q), \sharp^{1}_\mathscr{L}, - \can_{\sharp_\mathscr{L}}^{1}) ,\] 
		and deduce this result as a consequence of Lemma \ref{U0U1dual} applied to  Lemma \ref{lemma:symm-E}. 
	\end{proof} 
	\subsection{Twisted duality}  
	Let $\ell = 2a+1$ and let $\mathscr{L}= \mathscr{O}(-\ell)$.
	Define the functor 
	\begin{align*}
		\veet : ( \mathscr{O}\hModh C_0(q))^\op &\to \mathscr{O} \hModh C_0(q) \\
		M &\mapsto M^{\veet} 
	\end{align*}
	where $ M^{\veet} :=  [ C_1(q) \otimes_{C_0(q)} M^\op , \mathscr{L} ]_{\mathscr{O}}$. The double dual identification $\can:M \mapsto M ^{\flat \flat}$ is given by the formula
	\[ \can(x)(y \otimes g^\op) = (-1)^{|x||g|} \sigma(g(y \otimes x^\op))\]
	for $x \in M$, $y \in C_1(q)$ and $g \in [C_1(q) \otimes_{C_0(q)} M^\op, \mathscr{L}]_{\mathscr{O}}$. Construct an isomorphism 
	\begin{equation}\label{eq:twistedclifford}
		\xymatrix{
			\cdots \ar[r] & \mathscr{O}(-i-2) \otimes C_{i+r+2}(q) \ar[r] \ar[d]^{\theta}   & \mathscr{O}(-i-1)\otimes C_{i+r+1}(q) \ar[d]^{\theta} \ar[r] \ar[r] & \cdots \\
			\cdots \ar[r] & (\mathscr{O}(i+2-\ell) \otimes C_{i+r+1}(q))^{\veet} \ar[r] & (\mathscr{O}(i+1-\ell) \otimes C_{i+r}(q))^{\veet} \ar[r]  \ar[r] & \cdots
	}\end{equation}
	of chain complexes of $\mathscr{O}$-modules and right $C_0(q)$-modules 	where we define the map
	\[\theta: \mathscr{O}(-i) \otimes C_{r+i}(q) \longrightarrow (\mathscr{O}(i-\ell) \otimes C_{r+i-\ell}(q))^{\veet}  \]
	as 
	\[\theta (f \otimes x) ( c\otimes (g \otimes y)  ^\op) = fg ~ \tr(\sigma(x) \cdot  y \cdot  c).\]
	
	\begin{proposition}\label{prop:V0phi}
		Let $\ell = 2a+1$ . There is a non-degenerate symmetric form
		$$ \theta: \mathscr{U}_{a} \rightarrow (\mathscr{U}_{a} )^{\veet} $$ 
		in the category with duality $(\mathscr{O}\hModh C_0(q), \veet)$. 
	\end{proposition}

	\begin{proof}
		
		By taking cokernels of horizontal maps in \eqref{eq:twistedclifford}, we get isomorphisms
		\[\theta_i: \mathscr{U}_{i} \to (\mathscr{U}_{\ell-i-1})^\veet\]	
		and the form in the statement is defined by taking $i =a$.
		It is enough to check that the symmetry relation $ \theta^{\veet} \circ \varpi = \theta$. This follows from the symmetry of the trace map, i.e. $\tr(\sigma(x) \cdot  y \cdot  c)= \tr(\sigma(y) \cdot  x \cdot  c)$. 
	\end{proof}
	\begin{remark}\label{rmk:symm-C1}
		Alternatively, one may consider $\theta$ as a map
		\[ \theta: \mathscr{U}_{a} \otimes_{C_0(q)} C_1(q) \otimes_{C_0(q)} \mathscr{U}_{a}^\op \to \mathscr{L}  \]
		which is a symmetric form in the category $\mathscr{O}\hModh C_0(q)$ with respect to the duality coefficient in $(\mathscr{L}, \id)$ and $(C_1(q), \sigma)$ in sense of \cite[Section 7.6]{S-Inven}.
	\end{remark}
	
	\section{Proof of Theorem \ref{theo:main-theorem}}
	
	\subsection{Proof of the main theorem} The idea of the proof of Theorem \ref{theo:main-theorem} is to mutate the semi-orthogonal decomposition of $D^bQ_d$ into a form which respects the duality $\sharp_\mathscr{L}$, and then use the additivity theorem (cf. \cite[Lemma A.2]{KSW21}) and the localization theorem (cf. \cite[Theorem 6.6]{S-JPAA}) for Hermitian $K$-theory. We break down the proof Theorem \ref{theo:main-theorem} into the following pieces:
	\begin{enumerate}
		\item [(1)] 	Theorem \ref{theo:main-theorem} (i) will be proved in Theorem \ref{thm:(i)}.
		\item [(2)]	Theorem \ref{theo:main-theorem} (ii) will be proved in Theorem \ref{thm:(ii)}. 
		\item [(3)] Theorem \ref{theo:main-theorem} (a) will be proved in Proposition \ref{thm:C0C1S0}.
		\item [(4)] Theorem \ref{theo:main-theorem} (b) will be proved in Theorem \ref{theo:GW(S)}. 
	\end{enumerate}
	
	\subsection{Semi-orthogonal decomposition and duality}
	\begin{theo}
		There is a full exceptional collection in the sense of \cite{bondal-kapranov}
		\[D^bQ_d = \big\langle \mathscr{O}(1-d), \cdots, \mathscr{O}(-1-a), \mathscr{U}_a, \mathscr{O}(-a), \cdots,  \mathscr{O} \big\rangle\]
		for any integer $a$ such that $0\leq a < d-1$. 
	\end{theo}
	\begin{proof}
		The result follows from the proof of \cite[Corollary 3.8]{Xie19}.
	\end{proof}
	\begin{remark}\label{theo:sodofq}
		In this paper, we will use the following semi-orthogonal decompositions:
		\begin{itemize}
			\item [(i).] $ D^bQ=	\begin{cases}
				\langle	\mathcal{A}_{[1-d,-1-m]}, \mathcal{U}_{m}, \mathcal{A}_{-m}, \mathcal{A}_{[1-m,0]} \rangle  (\sharp_{\mathscr{O}(1-d)}) & \textnormal{if $d=2m+1$;}\\
				\langle  \mathcal{A}_{[1-d,-1-m]}, \mathcal{U}_{m},\mathcal{A}_{-m}, \mathcal{A}_{[1-m,-1]},  \mathcal{A}_{0} \rangle (\sharp_{\mathscr{O}(-d)})  & \textnormal{if $d=2m$.}
			\end{cases} $
			\item [(ii).]
			$ D^bQ=	\begin{cases}
				\langle	\mathcal{A}_{[1-d,-1-m]}, \mathcal{U}_{m}, \mathcal{A}_{[-m,0]} \rangle (\sharp_{\mathscr{O}(1-d)}) & \textnormal{if $d=2m+2$;}\\
				\langle \mathcal{A}_{[1-d,-1-m]}, \mathcal{U}_{m}, \mathcal{A}_{[-m,-1]},\mathcal{A}_0  \rangle (\sharp_{\mathscr{O}(-d)}) & \textnormal{if $d=2m+1$.}
			\end{cases} $
		\end{itemize}
		where the relevant duality is indicated in the round brackes (See \cite[Section 3]{Xie19} for the notations $\mathcal{A}_i, \mathcal{A}_{[i,j]}$ and $\mathcal{U}_i$). In \cite[Corollary 3.10]{Xie19}, we state the semi-orthogonal decomposition in a different form, which is because we implicitly use the ``two-periodicity in twists'' isomorphism in the proof of \cite[Theorem 1.2]{Xie19}. However, we will skip this extra step in this paper.
	\end{remark}
	
	\begin{definition}
		If $\mathscr{L} = \mathscr{O}(-\ell)$, we define the full triangulated subcategories $\mathcal{A}_{\mathscr{L}}$ of $ D^bQ$ as follows:
		\[\mathcal{A}_{\mathscr{L}}:= \begin{cases}
			\langle \mathscr{U}_{a} \rangle & \textnormal{if $[\mathscr{L}] = [\mathscr{O}(1)] $} \\
			\langle \mathscr{U}_{a}, \mathscr{O}(-a) \rangle & \textnormal{if $[\mathscr{L}] = [\mathscr{O}] $} 
		\end{cases}\]
		where $a = \lfloor \frac{\ell}{2} \rfloor$. 
	\end{definition}	
	
	\begin{lemma}\label{lem:A-fixed-duality}
		The  category $\mathcal{A}_\mathscr{L}$ is fixed by the duality $\sharp_\mathscr{L}$. 
	\end{lemma}		
	\begin{proof}
		If $\ell = 2a$, then
		\[
		\mathcal{A}_{\mathscr{L}} = \big\langle \mathscr{U}_{a}, \mathscr{O}(-a) \big\rangle .
		\]
		Clearly, $\mathscr{O}(-a)^{\sharp_{\mathscr{L}}} \cong \mathscr{O}(-a)$. The claim follows from the exact sequence
		\[
		\xymatrix{0 \ar[r] & \mathscr{U}_{a} \ar[r] & \mathscr{O}(-a)\otimes C_{r+a}(q) \ar[r] & (\mathscr{U}_a)^{\sharp_{\mathscr{L}}} \ar[r] & 0 }
		\]
		by noting that $(\mathscr{U}_a)^{\sharp_{\mathscr{L}}} \in \mathcal{A}_{\mathscr{L}}$. If $\ell = 2a+1$, then
		\[
		\mathcal{A}_{\mathscr{L}} = \big\langle \mathscr{U}_{a} \big\rangle .
		\]
		Now we see the claim since $(\mathscr{U}_{a})^{\sharp_{\mathscr{L}}}$ is isomorphic to $ \mathscr{U}_{a} \otimes_{C_0(q)} C_1(q)$ as  $\mathscr{O}$-modules by Proposition \ref{prop:V0phi}.
	\end{proof}
	
	\begin{definition}
		Let $\mathcal{A}$ be any full triangulated subcategory of $ D^bQ$. The full dg subcategory of $\mathrm{Ch}^bQ$ of objects the same as $\mathcal{A}$ which is a pretriangulated dg category will be denoted by a different font $\mathfrak{A}$ for simplicity. 
	\end{definition}		
	
	The dg category $\mathfrak{A}$ is a kind of dg enhancement of $\mathcal{A}$. 
	
	\begin{theo}\label{thm:(i)}
		If $\mathscr{L} = \mathscr{O}(1-d)$, then there is a stable equivalence of spectra
		$$  \xymatrix{ K(k)^{\lfloor \frac{d}{2} \rfloor } \oplus \GW^{[i]}(\mathfrak{A}_{\mathscr{L}}) \ar[r]^-{\simeq} &  \GW^{[i]}(Q_d,\mathscr{L}). }  $$
	\end{theo}
	\begin{proof}
		The proof follows from the additivity theorem in Hermitian $K$-theory. One needs to check that all assumptions are satisfied. Let us work with the semi-orthogonal decomposition
		\[
		D^bQ_d = \big\langle \mathcal{A}_{[1-d,-\lceil \frac{d}{2} \rceil]}, \mathcal{A}_{\mathscr{L}}, \mathcal{A}_{[1-\lfloor \frac{d}{2} \rfloor,0]} \big\rangle .
		\]
		By Lemma \ref{lem:A-fixed-duality}, the triangulated category $\mathcal{A}_{\mathscr{L}}$ is fixed by the duality $\sharp_{\mathscr{L}}$, and so is its dg enhancement $\mathfrak{A}_\mathscr{L}$.\ Since the categories $\mathcal{A}_{[1-\lfloor \frac{d}{2} \rfloor,0]}$ and $\mathcal{A}_{[1-d,-\lceil \frac{d}{2} \rceil]}$ are switched mutually by the duality $\sharp_{\mathscr{L}}$, the additivity theorem applies, and we deduce a stable equivalence
		\[
		K(\mathfrak{A}_{[1-\lfloor \frac{d}{2} \rfloor,0]}) \oplus \GW^{[i]}(\mathfrak{A}_{\mathscr{L}})
		\ \xrightarrow{\ \cong\ }\ \GW^{[i]}(Q_d,\mathscr{L})
		\]
		of spectra. Finally, we conclude by noting that
		\[
		K(\mathfrak{A}_{[1-\lfloor \frac{d}{2} \rfloor,0]})
		\cong \bigoplus_{i=1}^{\lfloor \frac{d}{2} \rfloor} K(\mathfrak{A}_{1-i})
		\cong K(k)^{\lfloor \frac{d}{2} \rfloor},
		\]
		where the first equivalence is obtained by the additivity theorem, and the second equivalence follows from Thomason's invariance theorem for algebraic $K$-theory (cf.\ \cite[Theorem 3.2.24]{S3} and \cite[Theorem 3.2.29]{S3}).
	\end{proof}
	
	\begin{theo}\label{thm:(ii)}
		If $\mathscr{L} = \mathscr{O}(-d)$, then there is a homotopy fibration sequence of spectra
		$$  \xymatrix{ K(k)^{\lceil \frac{d}{2} \rceil -1} \oplus \GW^{[i]}(\mathfrak{A}_{\mathscr{L}})  \ar[r] &   \GW^{[i]}(Q_d,\mathscr{L})  \ar[r]^-{p_*} &  \GW^{[i-d]}(k, \Lambda^n P). 	}	$$
	\end{theo}
	\begin{proof}
		The proof will use the localization theorem and the additivity theorem. Let us work with the semi-orthogonal decomposition
		\[D^bQ_d = \big\langle  \mathcal{B},  \mathcal{A}_{0} \big\rangle\]
		where 
		\[ \mathcal{B} = \big \langle \mathcal{A}_{[1-d,-1-\lfloor \frac{d}{2} \rfloor]}, \mathcal{A}_{\mathscr{L}}, \mathcal{A}_{[1-\lceil \frac{d}{2}\rceil,-1]} \big \rangle\]
		
		By Lemma \ref{lem:A-fixed-duality}, the category $\mathcal{B}$ is fixed by the duality $\sharp_{\mathscr{L}}$, and so is its dg enhancement $\mathfrak{B}$. 	Let $v$ be the set of morphisms in $\textnormal{Ch}^b Q_d$ that become isomorphisms in the Verdier quotient 
		$D^b Q_d  / \mathcal{B} . $  Therefore we obtain a quasi-exact sequence
		$$ \xymatrix{(\mathfrak{B},\quis, \sharp_{\mathscr{L}}) \ar[r] &  (\textnormal{Ch}^b Q_d, \quis, \sharp_{\mathscr{L}}) \ar[r]  & (\textnormal{Ch}^b Q_d, v, \sharp_{\mathscr{L}}) }$$
		of dg categories with duality, and the localization theorem yields a homotopy fibration
		\begin{equation}\label{WQACDbQ1} \xymatrix{\GW^{[i]}(\mathfrak{B}, \sharp_{\mathscr{L}}) \ar[r]  & \GW^{[i]}(Q_d, \mathscr{L}) \ar[r] & \GW^{[i]}(\textnormal{Ch}^b Q_d, v, \sharp_{\mathscr{L}}). } \end{equation}
		By the additivity theorem, we deduce stable equivalences 
		$$  \GW^{[i]}(\mathfrak{B}, \sharp_{\mathscr{L}})  \cong     K(\mathfrak{A}_{[1-\lceil \frac{d}{2}\rceil,-1]}) \oplus \GW^{[i]}(\mathfrak{A}_{\mathscr{L}}) \cong K(k)^{\lceil \frac{d}{2} \rceil-1} \oplus \GW^{[i]}(\mathfrak{A}_{\mathscr{L}}).  $$
		Finally, tensoring the form $(\mathbf{K},{\varphi })$ in Proposition \ref{prop:Kphi} yields a non-singular exact dg form functor
        \[(\mathbf{K},{\varphi }) \otimes-:(\mathrm{Ch}^bk, \quis, \sharp_{\Lambda^nP}) \to (\mathrm{Ch}^b Q_d, v, \sharp_{\mathscr{L}})  \]
        which induces a stable equivalence 
		$$ \xymatrix{ (\mathbf{K},{\varphi }) \otimes - : \GW^{[i-d]}(k,\, \Lambda^nP)  \ar[r]^-{\cong}  &  \GW^{[i]}(\textnormal{Ch}^b Q_d, v, \sharp_{\mathscr{L}}) }$$
		of spectra. To see this functor is non-singular, we note that the cone $C(\varphi)$ of $\varphi$ vanishes in $D^bQ_d/\mathcal{B}$ by Lemma \ref{lemma:symm-E}, since $\mathbf{E}$ is acyclic and $\mathbf{N}$ lands in $\mathcal{B}$.  
	\end{proof}
	\begin{remark}
		In the statement, we define the map $p_*$ by the commutativity of the following diagram
		\[ \xymatrix{ \GW^{[i]}(Q_d, \mathscr{L}) \ar[d]_-{p_*} \ar[r] & \GW^{[i]}(\textnormal{Ch}^b Q_d, v, \sharp_{\mathscr{L}})  \\
			\GW^{[i-d]}(k,\, \Lambda^nP) \ar[ur]^-{\cong}_-{\quad (\mathbf{K},{\varphi })\otimes- } }  \] 
		On the level of derived categories, the commutativity of the diagram
		\[ \xymatrix{D^bQ_d \ar[r] \ar[d]_-{Rp_*} & D^bQ_d / \mathcal{B} \\
			D^b k \ar[ur]_-{\mathbf{K}\otimes-} } \]
		of triangulated categories justifies the notation $p_*$, where we note that $D^bQ_d / \mathcal{B} \cong \mathcal{A}_0$ and $Rp_*(-) = R\Hom(\mathscr{O},-)$, and that $\mathbf{K}\cong \mathscr{O}$ in $D^bQ_d / \mathcal{B}$. The map $p_*$ coincides with the pushforward map constructed in \cite[Section 3]{HX23} if $k$ is regular. 
	\end{remark}
	
	\subsection{The case $[\mathscr{L}] = [\mathscr{O}(1)]$} Let $\ell = 2a+1$ and $\mathscr{L}:= \mathscr{O}(-\ell)$. By Proposition \ref{prop:V0phi}, the dg  category $\mathfrak{A}_{\mathscr{L}}$ is fixed by the duality $\sharp_\mathscr{L}$. By our convention, the dg category $\mathfrak{A}_{\mathscr{L}}$ has the same objects as the triangulated category $\langle \mathscr{U}_a \rangle$. 
	\begin{proposition}\label{thm:C0C1S0} The map of spectra 
		$$ \xymatrix{ (\mathscr{U}_{a} , \theta)\otimes-: \GW^{[i]}(C_0(q)_\sigma, C_1(q)_\sigma) \xrightarrow{\quad \cong \quad} \GW^{[i]}(\mathfrak{A}_{\mathscr{L}}, \sharp_{\mathscr{L}}) }$$
		is a stable equivalence in $\mathcal{SH}$.
	\end{proposition}
	\begin{proof}
		By \cite[Section 7.7, p.  404]{S-Inven} and Remark \ref{rmk:symm-C1}, we get a dg form functor
		\[  (\mathscr{U}_{a} , \theta)\otimes-: (C_0(q)\hMod, \sharp_{C_1(q)_\sigma}) \to (\mathscr{O}\hMod,\sharp_{\mathscr{L}}) \]
		which induces a non-singular exact dg form functor 
		\[  (\mathscr{U}_{a} , \theta)\otimes-: (\mathrm{Ch}^bC_0(q), \quis,\sharp_{C_1(q)_\sigma}) \to (\mathfrak{A}_{\mathscr{L}},\quis, \sharp_{\mathscr{L}}) \]
		by Proposition \ref{prop:V0phi}, since $\theta$ is an isomorphism. 
	\end{proof}
	\subsection{The case $[\mathscr{L}] = [\mathscr{O}]$} 
	Let $\ell = 2a$ and $\mathscr{L} = \mathscr{O}(-\ell)$. 		Recall the dg form functor 
	$$
	(k \hModh C_0(q), \sharp_k)\otimes (C_0(q) \hMod, \sharp_{C_0(q)_\sigma} ) \rightarrow (k \hMod , \sharp_k ) 
	$$
	cf. \cite[Lemma 2.19]{Xie19}. 
	Tensoring the symmetric form 
	$$\vartheta_{r+a}: C_{r+a}(q) \rightarrow \Hom_k(C_{r+a}(q)^\op, k), \quad x \mapsto (y \mapsto \tr(\sigma(x)\cdot y))$$ 
	in $(k \hModh C_{0}(q), \sharp_k)$ 
	induces a dg form functor
	$$\xymatrix{ (C_{r+a}(q), \vartheta_{r+a}) \otimes_{C_0(q)} - : ( C_0(q) \hMod, \sharp_{C_0(q)_\sigma}  )  \longrightarrow (k \hMod , \sharp_k  ) } $$ 
	cf. \cite[Lemma 2.17]{Xie19}, which gives rise to a map of Grothendieck-Witt spectrum 
	\begin{equation}\label{eq:traceC}
		\tr_\ell: \GW^{[i]}(C_0(q)_\sigma) \rightarrow \GW^{[i]}(k)\end{equation} 
	(Compare to \cite[Lemma 6.8]{Xie19} for the case of Witt groups). 		
	
	\begin{theo}\label{theo:GW(S)}
		Suppose that $[\mathscr{L}] = [\mathscr{O}]$. Then, there is a distinguished triangle
		\begin{equation}\label{eq:dist-triangle}
			\xymatrixcolsep{2pc} \xymatrix{  \GW^{[i]}(k)  \ar[r] & \GW^{[i]}(\mathfrak{A}_{\mathscr{L}}, \sharp_{\mathscr{L}})   \ar[r]^-{  } & \GW^{[i+1]}(C_0(q)_\sigma)  \ar[r]^-{  \eta \cup \tr_\ell } &   \GW^{[i]}(k) \big[ 1\big] }
		\end{equation} 
		in the stable homotopy category of spectra, where  $\eta \in \GW^{[-1]}_{-1}(k)$ is the Bott element defined in \cite[(6.1) p.\ 56]{S-JPAA}. 
	\end{theo}
	\begin{proof}
		
		Let $w$ be the set of morphisms in the dg category $\mathfrak{A}_{\mathscr{L}}$  that become isomorphisms in the Verdier quotient 
		$\mathcal{A}_{\mathscr{L}}  /\mathcal{A}_{-a} $.\ It follows that the sequence 
		$$ \xymatrix{   (\mathfrak{A}_{-a} , \quis, \sharp_\mathscr{L}) \longrightarrow (\mathfrak{A}_{\mathscr{L}}, \quis,\sharp_\mathscr{L}) \longrightarrow ( \mathfrak{A}_{\mathscr{L}}, w, \sharp_\mathscr{L}) }$$
		is quasi-exact, which gives rise to a homotopy fibration 
		\begin{equation}\label{GWloca}  \xymatrix{  \GW^{[i]}(\mathfrak{A}_{-a} , \sharp_\mathscr{L}) \ar[r] &  \GW^{[i]}(\mathfrak{A}_{\mathscr{L}},\sharp_\mathscr{L})  \ar[r] & \GW^{[i]}  ( \mathfrak{A}_{\mathscr{L}}, w, \sharp_\mathscr{L}) } \end{equation} 
		by the localization theorem \cite{S-JPAA}. 			
		Since the exact dg form functor 
		\[ \mathscr{O}(-a) \otimes -: (\textnormal{Ch}^b k,\quis, \sharp_k) \longrightarrow ( \mathfrak{A}_{-a}, \quis, \sharp_\mathscr{L})\]
		induces an equivalence of associated triangulated categories,  we deduce a stable equivalence of spectra
		\begin{equation}\label{eq:o1}
			\xymatrix{\mathscr{O}(-a) \otimes -: \GW^{[i]}(k) \xrightarrow{ \cong } \GW^{[i]}( \mathfrak{A}_{-a}, \sharp_\mathscr{L}) }		
		\end{equation}
		by the invariance theorem, cf.  \cite[Theorem 6.5]{S-JPAA}. 
		Moreover, Proposition \ref{prop:clmu} brings out a non-singular exact dg form functor
		\[ Cl \otimes - : (\mathrm{Ch}^bC_0(q), \quis,  \sharp_{C_0(q)_\sigma}^{1}) \to (\mathfrak{A}_\mathscr{L},w, \sharp_{\mathscr{L}}) \]
		which induces a stable equivalence
		\begin{equation}\label{GWCA} Cl\otimes- : \GW^{[i+1]}(C_0(q)_\sigma) \xrightarrow{ \cong } \GW^{[i]}(\mathfrak{A}_\mathscr{L},  w, \sharp_\mathscr{L})\end{equation}
		(Compare to \cite[Proposition 6.6]{Xie19}). Putting (\ref{eq:o1}) and (\ref{GWCA}) into the homotopy fibration (\ref{GWloca}), we obtain the distinguished triangle
		\[ \xymatrixcolsep{2pc} \xymatrix{  \GW^{[i]}(k)  \ar[r] & \GW^{[i]}(\mathfrak{A}_{\mathscr{L}}, \sharp_{\mathscr{L}})   \ar[r]^-{  } & \GW^{[i+1]}(C_0(q)_\sigma)  \ar[r]^-{  \partial_\mathscr{L} } &   \GW^{[i]}(k) \big[ 1\big] } \]
		It remains to show that the boundary map $\partial_\mathscr{L}$ is homotopic to $\eta \cup \tr_\ell$, which follows from Lemma \ref{boundarytracemap} below.   \end{proof}
	
	\begin{lemma}\label{boundarytracemap} The diagram 
		\[  \xymatrix{  \GW^{[i+1]}(C_0(q)) \ar[r]^-{\tr_\ell} \ar[d]_-{\partial_\mathscr{L}} & \GW^{[i+1]}(k) \ar[dl]^-{\quad \eta\cup} \\
			\GW^{[i]}(k) \big[1\big]  }\]
		commutes in the stable homotopy category of spectra. 
	\end{lemma}
	\begin{proof}
		Let $\Fun([1], \mathfrak{A}_\mathscr{L}) $ be the dg category consisting of objects $A = (A_0 \rightarrow A_1) \in \mathfrak{A}_\mathscr{L}$  and morphisms 
		$$ \xymatrix{ A\ar[d]^-{f} & A_0  \ar[r] \ar[d]^-{f_0}& A_1 \ar[d]^-{f_1}\\
			B& B_0 \ar[r] & B_1}$$ 
		(cf.\ Section 2.1 \cite{S-JPAA}).\  The unit map $[1] \rightarrow [0]$ induces an exact dg functor $I: \mathfrak{A}_\mathscr{L}\rightarrow \Fun([1], \mathfrak{A}_\mathscr{L}) \colon A \rightarrow 1_A$.\ A quasi-isomorphism $f:A\rightarrow B$ in $ \Fun([1], \mathfrak{A}_\mathscr{L}) $ means $f_0$ and $f_1$ are both quasi-isomorphisms of complexes.\ 
		Let $ \Fun_{w}([1], \mathfrak{A}_\mathscr{L}) $ be the dg full subcategory of $ \Fun([1], \mathfrak{A}_\mathscr{L}) $ consisting of objects $A_0\rightarrow A_1 \in w$.\ Define $P: \Fun_{w}([1], \mathfrak{A}_\mathscr{L})  \rightarrow \mathfrak{A}$ by sending $A_0\rightarrow A_1$ to $A_0 $.\ Consider the following diagram consisting of exact dg form functors
		\begin{equation}\label{BS-Inven} 
			\xymatrix{ (\mathfrak{A}^w_{\mathscr{L}} ,\quis, \sharp_\mathscr{L}) \ar[r]^-{I} \ar[d] & ( \Fun([1], \mathfrak{A}_\mathscr{L}^w),\quis, \sharp_\mathscr{L}) \ar[r]^-{\cone} \ar[d] &  (\mathfrak{A}^w_{\mathscr{L}} ,\quis, \sharp_\mathscr{L}^1) \ar^-{=}[d] \\ 
				(\mathfrak{A}_{\mathscr{L}},  \quis, \sharp_\mathscr{L}) \ar[r]^-{I} \ar[d] & ( \Fun_{w}([1], \mathfrak{A}_\mathscr{L}) , \quis, \sharp_\mathscr{L}) \ar[r]^-{\cone} \ar[d]^-P & (\mathfrak{A}^w_{\mathscr{L}} ,\quis, \sharp_\mathscr{L}^1) \\
				(\mathfrak{A}_{\mathscr{L}},  w, \sharp_\mathscr{L}) \ar[r]^-{=} &  (\mathfrak{A}_{\mathscr{L}}, w , \sharp_\mathscr{L}) }\end{equation} 
		where all the vertical and horizontal sequences are quasi-exact.\ Taking $\GW$-spectra, we obtain the following commutative diagram  
		\begin{equation}\label{BS-JPAA}  \xymatrix{ \GW^{[i]}(\mathfrak{A}^w_{\mathscr{L}} ) \ar[r]^-{I} \ar@{}[dr]|{\square} \ar[d] & \GW^{[i]}( \Fun([1], \mathfrak{A}_\mathscr{L}^w)) \ar[r]^-{\cone} \ar[d] &  \GW^{[i+1]}(\mathfrak{A}^w_{\mathscr{L}} ) \ar^-{=}[d]  \ar[r]^-{-\eta\cup} &  \GW^{[i]}(\mathfrak{A}^w_{\mathscr{L}} ) \big[1\big] \\ 
				\GW^{[i]}(\mathfrak{A}_{\mathscr{L}})  \ar[r]^-{I} \ar[d] & \GW^{[i]}( \Fun_{w}([1], \mathfrak{A}_\mathscr{L}) ) \ar[r]^-{\cone} \ar[d]^-P & \GW^{[i+1]}(\mathfrak{A}^w_{\mathscr{L}} ) \\
				\GW^{[i]}(\mathfrak{A}_{\mathscr{L}}, w ) \ar[r]^-{=} \ar[d]_-{\partial} &  \GW^{[i]}(\mathfrak{A}_{\mathscr{L}}, w ) \\
				\GW^{[i]} (\mathfrak{A}^w_{\mathscr{L}} ) \big[1\big] }\end{equation} 
		in the stable homotopy category  of spectra where the diagram $\square$ is a homotopy pushout.\ By \cite[Appendix A]{HX23}, the following diagram commutes 
		\begin{equation}\label{eq:Fun} \xymatrix{ \GW^{[i]}( \Fun_{w}([1], \mathfrak{A}_\mathscr{L}) ) \ar[r]^-{\cone} \ar[d]^-{P} & \GW^{[i+1]}(\mathfrak{A}^w_{\mathscr{L}} ) \ar[d]^-{\eta\cup} \\
				\GW^{[i]}(\mathfrak{A}_{\mathscr{L}}, w ) \ar[r]^-{\partial} &   \GW^{[i]} (\mathfrak{A}^w_{\mathscr{L}} ) \big[1\big]  }\end{equation}
		
		Diagram (\ref{eq:Fun}) fits into the following commutative diagrams. 
		\begin{equation}\label{eq:trace-fibration}
			\xymatrix{  \GW^{[i+1]}(C_0(q)_\sigma) \ar[d]^-{\mu\otimes - } \ar@/_5pc/[dd]_{Cl\otimes - }       \ar@{=}[rr] && \GW^{[i+1]}(C_0(q)_\sigma) \ar[d]^\tr \\
				\GW^{[i]}( \Fun_{w}([1], \mathfrak{A}_\mathscr{L}) ) \ar[r]^-{\cone} \ar[d]^-{P}  & \GW^{[i+1]} (\mathfrak{A}^w_{\mathscr{L}} )\ar[d]^-{\eta\cup} & \ar[l]_-{\cong}  \GW^{[i+1]}(k) \ar[d]^-{\eta\cup} \\
				\GW^{[i]}(\mathfrak{A}_{\mathscr{L}}, w ) \ar[r] &   \GW^{[i]} (\mathfrak{A}^w_{\mathscr{L}} ) \big[1\big] & \ar[l]_-{\cong}    \GW^{[i]} (k) \big[1\big]  }
		\end{equation} 
		where $\mu = (\mu, \Sigma) $  is defined by the diagram after \cite[Formula (26), page 100]{Xie19}. Note that $P \circ (\mu\otimes - ) = Cl \otimes -$ by the symmetry of $B$.\ Moreover, the upper square in Diagram \eqref{eq:trace-fibration}
		commutes up to homotopy. To see this, we need to construct a natural transformation 
		$$ F:  (\mathscr{O} \otimes-) \circ \tr \rightarrow \cone \circ (  \mu \otimes - )$$ of exact dg form functors with $F_M $ a quasi-isomorphism for any $M \in \mathrm{Ch}^bC_0(q)$, which is already done in the proof of \cite[Theorem 6.11, Step II]{Xie19}.
	\end{proof}

	\begin{appendix}
		
		\section{Truncating symmetric forms}\label{sec:symmform}
		For preliminaries on (exact) categories with dualities and sign choices, we refer to \cite[Section 6]{S-Inven}. Let $(\mathcal{E},* , \eta)$ be an exact category with duality. Recall from \cite[Section 6.1]{S-Inven} that the duality $(*,\eta)$ induces the (naive) duality $(*^n,\eta^n)$ on $\mathrm{Ch}^b\mathcal{E}$ by the formulas in \cite[Section 6.1 p. 393]{S-Inven}. This naive duality will be convenient when we wish to construct explicit symmetric forms on $\mathrm{Ch}^b\mathcal{E}$. \footnote{There is another equivalent duality on $\mathrm{Ch}^b \mathcal{E}$, namely $\sharp^n$, defined in \cite[Section 6.1, Remark 13, p. 394]{S-Inven}. In this paper, we use both dualities. The duality $*^n$ is only used in Section \ref{sub:tate-dual}. } The following result plays an important role in constructing symmetric forms. 
		\begin{lemma}\label{lemma:symm-E}
			Let $(\mathscr{E},*,\eta)$ be an exact category with duality, and let $\epsilon= \pm 1 $. Let 
			$$ \xymatrix{ \mathbf{E} \ar[d]_{\phi} &  E^{-n} \ar[d]_-{\phi_{-n}} \ar[r]^-{d^{-n}}  &  E^{-n+1} \ar[d]_-{\phi_{-n+1}} \ar[r] &\cdots \ar[r] & E^{0} \ar[d]_-{\phi_{0}} \ar[r]^-{d^0}  & \ar[d]_-{\phi_{1}} E^{1} \\
				\quad	\mathbf{E}^{*^{n-1}} &  (E^{1})^* \ar[r]^-{(d^{0})^*}  &  (E^{0})^* \ar[r] &   \cdots \ar[r] & (E^{-n+1})^* \ar[r]^-{(d^{-n})^*}  & (E^{-n})^*  }  $$ 
			be an $\epsilon$-symmetric form in the category with duality $(\mathrm{Ch}^b(\mathscr{E}), *^{n-1}, \eta^{n-1}  )$ such that $\phi_1$ is an isomorphism.  
			Then, truncating $E^1$ yields a $(-\epsilon)$-symmetric form 
			$$ \xymatrix{
				\mathbf{K} \ar[d]_-{{\varphi }} &  E^{-n} \ar[d]_-{\varphi_{-n}} \ar[r]^-{d^{-n}}  &  E^{-n+1} \ar[d]_-{\varphi_{-n+1}}  \ar[r] &   \cdots \ar[r] &  E^{-1} \ar[d]_-{\varphi_{-1}} \ar[r]^-{d^{-1}}  & E^{0} \ar[d]_-{\varphi_{0}} \\
				\quad	\mathbf{K}^{*^{n}} &  (E^{0})^{*} \ar[r]^-{(d^{-1})^{*}}  &  (E^{-1})^{*} \ar[r] &   \cdots \ar[r] &  (E^{-n+1})^{*} \ar[r]^-{(d^{-n+1})^{*}}  & (E^{-n})^{*} }$$ 
			in the category with duality $(\mathrm{Ch}^b(\mathscr{E}), *^{n}, \eta^{n})$ if  we let  
			$$\varphi_0= (-1)^n(d^{-n})^*  \phi_{0}/2, \quad  \varphi_{-n} = (d^{0})^* \phi_{-n}/2 \textnormal{ \quad and \quad }  \varphi_{i} = 0 \textnormal{ for $i \neq 0, -n$.}$$  
			Moreover, the cone of ${\varphi }$ fits into an exact sequence of chain complexes
			$$\xymatrix{\mathbf{E}[1] \ar[r] &  {C}({\varphi }) \ar[r] & \mathbf{N}}$$
			where $\mathbf{N}$ is the truncated chain complex 
			$$\xymatrix{ \mathbf{N} &       (E^{0})^* \ar[r]^-{(d^{-1})^*} & (E^{-1})^* \ar[r] & \cdots \cdots \ar[r] & (E^{-n+2})^* \ar[r]^-{(d^{-n})^*} & (E^{-n+1})^* }$$
			with $(E^{-n+1})^*$ in degree $-1$. 
		\end{lemma}
		
		\begin{proof}
			
			Note that $(d^{**})\eta^{n-1}_E = (\eta^{n-1}_E) d$, $ d^{*^{n-1}} \phi = \phi d$ and $\epsilon \phi^{*^{n-1}} \eta^{n-1} =\phi .$
			Note that $\eta^{n-1}$ is a map of chain complexes, so we have $$ (d^{-n})^{**}(\eta^{n-1}_E)_{-n}  = (d^{*^{n-1}*^{n-1}})_{-n}(\eta^{n-1}_E)_{-n} = (\eta^{n-1}_E)_{-n+1} d_{-n}.$$
			We show that $\varphi$ is $(-\epsilon)$-symmetric in the category with duality $(\mathrm{Ch}^b(\mathscr{E}), *^{n}, \eta^{n})$. This follows from the following identities.
			$$\begin{array}{lll} 
				\varphi_{-n} & = & \phi_{-n+1}d^{-n} /2
				\\
				& = &   \epsilon (\phi^{*^{n-1}} )_{-n+1}(\eta^{n-1}_E)_{-n+1} d^{-n} /2
				\\
				& = &   \epsilon (\phi_{0}																											  )^*(\eta^{n-1}_E)_{-n+1} d^{-n} /2
				\\ 
				& = &   \epsilon (\phi_{0}																											  )^* (d^{-n})^{**}(\eta^{n-1}_E)_{-n} /2
				\\
				& = &   \epsilon (-1)^{\frac{(n-1)(n-2)}{2}} (\phi_{0}																											  )^* (d^{-n})^{**}(\eta_{E_{-n}})/2
				\\ 
				& = &   -\epsilon (-1)^{\frac{n(n-1)}{2}}						((-1)^{n}(d^{-n})^{*}\phi_{0})^{*}(\eta_{E_{-n}})/2
				\\
				& = & -\epsilon (-1)^{\frac{n(n-1)}{2}} (\varphi_{0})^{*} (\eta_{E_{-n}})
				\\
				& = & - \epsilon (\varphi^{*^{n}})_{-n} (\eta^{n}_E)_{-n}
			\end{array} $$
			Similar symmetry holds for $\varphi_0$. 
			%
			For the second claim, we construct an exact sequence of chain complexes 
			$$\xymatrix{\mathbf{E}[1] \ar[r] &  \mathbf{C}({\varphi }) \ar[r] & \mathbf{N}}$$
			by the ladder diagrams
			$$ \xymatrix@C=1.5pc{
				E^{-n} \ar[d]^-{\id} \ar[r]^-{-d^{-n}}  &  E^{-n+1} \ar[d]^-{\kappa_{-n}}  \ar[r]^-{-d^{-n+1}} & E^{-n+2}  \ar[r] \ar[d]^-{\kappa_{-n+1}} & \cdots \ar[r] &  E^{0} \ar[r]^-{-d^0} \ar[d]^-{\kappa_{-1}} & E^{1} \ar[d]_-{(-1)^{n+1}\phi_1} 
				\\
				E^{-n} \ar[r]^-{\tilde{\Delta}} \ar[d] & E^{-n+1} \oplus (E^{0})^{*}  \ar[r]^{\Delta_{-n}}\ar[d]^-{\wp_{-n}} &  E^{-n+2} \oplus (E^{-1})^{*} \ar[d]^-{\wp_{-n+1}} \ar[r] &  \cdots \ar[r]^-{\Delta_{-2}} & E^0 \oplus  (E^{-n+1})^{*}  \ar[d]^-{\wp_{-1}} \ar[r]^-{\Delta}  & (E^{-n})^{*}\ar[d]
				\\
				0 \ar[r] & (E^{0})^{*}  \ar[r]^-{(d^{-1})^*}  &  (E^{-1})^* \ar[r]^-{(d^{-2})^*} &  \cdots \ar[r] &  (E^{-n+1})^{*}  \ar[r]  & 0
			}$$
			where we define
			$$\Delta_{-i} = \begin{pmatrix}
				-d^{-i+1} & 0 \\ 
				0 &  (d^{i-n-1})^*
			\end{pmatrix},\quad \Delta = \begin{pmatrix}
				(-1)^n\varphi_0 & (d_{-n})^* 
			\end{pmatrix}, \quad \tilde{\Delta} = \begin{pmatrix}
				-d^{-n} \\ \varphi_{-n} 
			\end{pmatrix}, $$
			$$\kappa_{-i} = \begin{pmatrix}
				1 \\ (-1)^{n-i-1} \phi_{-i+1}/2
			\end{pmatrix},\quad \wp_{-i} = \begin{pmatrix}
				(-1)^{n-i} \phi_{-i+1}/2 & 1 
			\end{pmatrix}$$
			This finishes the proof. 
		\end{proof}
	\end{appendix}

	\begin{ackno} Theorem \ref{theo:main-theorem} for the case of $[\mathscr{L}] = [\mathscr{O}]$ forms part of my PhD thesis in the University of Warwick. It is my pleasure to thank my supervisor Marco Schlichting for his insights. The author is partially supported by National Key Research and Development Program of China No.\ 2023YFA1009801, NSFC Grant 12271529, NSFC Grant 12271500.   The author  would also like to acknowledge the EPSRC
		Grant EP/M001113/1, DFG Priority Programme SPP 1786 and the DFG-funded research training group
		GRK 2240: Algebro-Geometric Methods in Algebra, Arithmetic and Topology.
	\end{ackno}


\begin{thebibliography}{9}
		
		\bibitem{bondal-kapranov}
		A. Bondal and M. Kapranov, 
		{\it Representable functors, Serre functors, and reconstructions,}
		Mathematics of the USSR Izvestia 35 (1990), 519–541.
		
		\bibitem{HX23}
		T. Huang and H. Xie, 
		{\it  The connecting homomorphism in Hermitian $K$-theory\/},
		preprint (2024).

        		
		\bibitem{KSW21}
		M. Karoubi, M. Schlichting and C. Weibel, 
		{\it  Grothendieck-Witt groups of some singular schemes\/},
       Proc. Lond. Math. Soc. (4) 122  (2021), 521-536.   
       
		\bibitem{S-Inven}
		M. Schlichting, 
		{\it  The Mayer-Vietoris principle for Grothendieck-Witt groups of schemes\/},
		Invent. Math. 179 (2010), no. 2, 349-433.
		\bibitem{S3}
		M. Schlichting, 
		{\it  Higher Algebraic $K$-theory (After Quillen, Thomason and Others)\/},
		Springer Lecture Notes in Math. 2008 (2011), 167-242.
		\bibitem{S-JPAA}
		M. Schlichting, 
		{\it Hermitian $K$-theory, derived equivalences, and Karoubi's fundamental theorem, \/} J. Pure Appl. Algebra (special issue in honor of Charles A. Weibel's 65th birthday).
		
		\bibitem{Swan}
		R. G. Swan,
		{\it $K$-theory of quadric hypersurfaces,\/}
		Ann. of Math., 122 (1985), 113-153.
		
		\bibitem{Xie19}
		H. Xie,
		{\it Witt groups of smooth projective quadrics,\/}
		Advances in Mathematics 346 (2019), 70-123. 
		\bibitem{Zibrowius}
		M. Zibrowius,  
		{\it Witt groups of complex cellular varieties,\/}
		Documenta Mathematica 16 (2011), 465–511. 
	\end{thebibliography}
\end{document}